\documentclass[final,3p,times]{elsarticle}
\usepackage[utf8]{inputenc}
\usepackage[hidelinks]{hyperref}

\usepackage{natbib}
\usepackage{graphicx}
\usepackage{amsmath}
\usepackage{amsfonts}
\usepackage{amsthm}
\usepackage{todonotes}
\usepackage{multirow}
\usepackage{layouts}
\usepackage{subcaption}
\usepackage{pgfplots}
\usetikzlibrary{external}
\usepgfplotslibrary{external}%
\newcommand{\R}{\mathbb{R}}
\newcommand{\Lkal}{\mathcal{L}}
\newcommand{\NN}{\textnormal{NN}}
\newcommand{\train}{\textnormal{train}}
\DeclareMathOperator{\argmin}{argmin}

\usepackage{mathtools}
\usepackage{etoolbox}

\DeclarePairedDelimiterX\set[2]{\{}{\}}
  {#1 \mathrel{}\mathclose{}\delimsize|\mathopen{}\mathrel{} #2}

\newtheorem{thm}{Theorem}
\newtheorem{lem}[thm]{Lemma}

\def\statement{\begin{minipage}[t]{.75\textwidth}
       Preprint
       \end{minipage}}

\makeatletter
\def\ps@pprintTitle{%
     \let\@oddhead\@empty
     \let\@evenhead\@empty
     \def\@oddfoot{\footnotesize\itshape
       \statement\hfill\today}%
     \let\@evenfoot\@oddfoot}
\makeatother

\begin{document}

\begin{frontmatter}

\title{Preconditioned FEM-based Neural Networks for Solving Incompressible Fluid Flows and Related Inverse Problems}
\author[label1]{Franziska Griese}
 \ead{franziska.griese@dlr.de}
 \cortext[cor1]{Corresponding author.}
\author[label1]{Fabian Hoppe}
 \ead{fabian.hoppe@dlr.de}
\author[label1]{Alexander Rüttgers}
 \ead{alexander.ruettgers@dlr.de}
\author[label1]{Philipp Knechtges\corref{cor1}}
 \ead{philipp.knechtges@dlr.de}
 \affiliation[label1]{organization={German Aerospace Center (DLR), Institute of Software Technology, High-Performance Computing Department},%
             addressline={Linder~Höhe}, 
             city={Cologne},
             postcode={51147}, 
             country={Germany}}

\begin{abstract}
The numerical simulation and optimization of technical systems described by partial differential equations is expensive, especially in multi-query scenarios in which the underlying equations have to be solved for different parameters. A comparatively new approach in this context is to combine the good approximation properties of neural networks (for parameter dependence) with the classical finite element method (for discretization). However, instead of considering the solution mapping of the PDE from the parameter space into the FEM-discretized solution space as a purely data-driven regression problem, so-called physically informed regression problems have proven to be useful. In these, the equation residual is minimized during the training of the neural network, i.e.\ the neural network "learns" the physics underlying the problem. In this paper, we extend this approach to saddle-point and non-linear fluid dynamics problems, respectively, namely stationary Stokes and stationary Navier--Stokes equations. In particular, we propose a modification of the existing approach: Instead of minimizing the plain vanilla equation residual during training, we minimize the equation residual modified by a preconditioner. By analogy with the linear case, this also improves the condition in the present non-linear case. Our numerical examples demonstrate that this approach significantly reduces the training effort and greatly increases accuracy and generalizability. Finally, we show the application of the resulting parameterized model to a related inverse problem.
\end{abstract}


\end{frontmatter}

\section{Introduction}

Understanding physical systems is crucial in many applications in science and engineering. Many important physical systems are modeled by partial differential equations (PDEs) and then solved using numerical methods. A more recent approach to simulate physical systems is to use neural networks (NNs) as surrogate models, which reduce complexity and speed up evaluation compared to numerical methods.
These deep neural networks are universal function approximators \cite{hornik1989multilayer} and are typically trained with data from simulations or measurements. However, generating data from experiments or simulations is expensive. In addition, a data-driven approach does not or only insufficiently take into account natural laws such as the conservation of energy, mass and momentum, thus generalizes poorly to unseen cases. To cope with these and similar issues, the field of physics-informed machine learning has emerged at the intersection of (numerical) mathematics, machine learning/artificial intelligence, and computational engineering; see, e.g., the surveys \cite{Karniadakis2021,Cuomo2022,Hao2023,Tanyu2023}.  

Among the numerous approaches to utilize the power of neural networks for solving differential equations, in particular so-called physics-informed neural networks (PINNs), first introduced in \cite{Dissanayake1994,Lagaris1998} and brought to broad attention by \cite{RAISSI2019686}, gained a lot of attention in the last years. PINNs embed the residual of the respective differential equations into the loss function of a neural network and compute the required derivatives exactly by automatic differentiation \cite{griewank2014automatic}. Hence, this approach is often considered to be able to avoid direct discretization; in particular, it is clearly a mesh-free approach. Since training exploits the underlying physics, encoded in the loss function, the PINN approach is able to solve a differential equation without the need for simulated or measured data. However, as higher order derivatives of the network in the loss function lead to a high learning complexity, and additional loss terms that incorporate boundary conditions lead to a multi-objective optimization problem (at least in the vanilla PINN approach), training PINNs is known to be potentially challenging; see, e.g., \cite{Wang2021}. For further literature dealing with PINNs, extensions thereof, or similar techniques such as the Deep Ritz method, the Deep Nitsche method, and weak adversarial networks --- all of these methods have in common that they do not explicitly apply discretization --- we refer to the aforementioned survey papers. 

In the following, we focus our literature overview on physics-informed neural networks that make explicit use of discretization, at least to a certain degree. The so-called variational PINNs (vPINNs), a Petrov--Galerkin-like approach, in which a neural network serves as trial function that is tested with a set of finitely many test functions, are introduced in \cite{kharazmi2019variational}. In the VarNet approach \cite{khodayi2020varnet}, specifically finite element functions are used as test functions, and in $hp$-vPINNs \cite{Kharazmi2021} essentially $hp$-finite element functions are used as test functions. In contrast, the Galerkin neural network approach \cite{ainsworth2021galerkin} extends the classical Galerkin procedure, that is usually bound to trial and test functions coming from linear spaces, to trial and test functions both coming from non-linear sets of functions parameterized by neural networks.
PhyGeoNet \cite{Gao2021} maps the underlying problem to a rectangular reference domain and then employs PINNs loss evaluated on a rectangular grid of points together with convolutional neural networks. Extending this idea to less structured meshes, physics-informed graph neural Galerkin networks \cite{gao2022physics} combine graph convolutional NNs with a Galerkin finite element method (FEM) discretization. In \cite{uriarte2022finite}, finite element discretization is considered and the proposed network architecture mimics subsequent refinement steps of the underlying mesh; moreover, a block Jacobi preconditioner is applied. In \cite{Songming2024} a somehow hybrid combination of PINNs and FEM-discretization is combined with preconditioning and applied to several model problems, including stationary Navier--Stokes equations (with homogeneous boundary conditions): while the network still approximates the non-discretized solution, it is only evaluated at grid points of a Lagrange FEM-discretization, and the loss functions for the training is constituted in terms of the left-preconditioned FEM-residual. Nonlinear problems are treated by iteratively solving linear problems. 
A model problem very similar to the one under consideration in the present paper is dealt with in \cite{Cao2025} using classical PINNs again. A so-called time-stepping-oriented training procedure, inspired by an control-theoretic ansatz to improve condition number of Jacobi matrices, is employed and allows to deal with 3d problems up to Reynolds number 5,000. We see these results as an encouraging sign that investigating the use of classical numerical techniques and ``tricks" is promising in the context of physics-informed ML. 

Recently, the use of PINNs has also been proposed to speed up certain element-local computations within a finite element method \cite{Pantidis2023}, and neural networks have been used to correct coarse level FEM-solutions on a finer level; see, e.g.,  \cite{Margenberg2022,Kapustsin2023,hintermueller2024}. Finally, a combination of Isogeometric Analysis with classical PINNs has been proposed in \cite{Moeller2021}. At the end of this overview, let us briefly comment also on two methods with a similar title that are only slightly related to the present context. In \cite{Mitusch2021}, the differential equation is discretized using classical FEM, while a neural network is used to parameterize an unknown parameter function in the differential equation. In \cite{Ramuhalli2005}, the finite element method is cast into a neural network. The domain-decomposition-based (nonlinear) preconditioning of the training of physics-informed neural networks developed in \cite{Kopanivcakova2024} refers to layer-parallel decomposition of the respective neural networks. 

Finally, \cite{meethal2023finite,LeDuc2023} combine neural networks and finite elements directly: the parametric differential equation under consideration is discretized a priori by finite elements, and a neural network is trained to map the parameters (of the differential equation) to the coefficients of the corresponding discrete solution (w.r.t.\ the chosen finite element basis). The loss function is given by the residual of the FEM-discretized differential equation, making this approach data-free and physics-informed. In \cite{meethal2023finite} the applicability to 1D-convection-diffusion problems, as well as to a truss problem, and in \cite{LeDuc2023} to further elasticity problems, such as the Euler--Bernoulli beam or the Reissner--Mindlin plate, has been demonstrated. For a discussion of the advantages of this approach compared to PINNs, in particular the simplified handling of boundary conditions or irregularly shaped domains, we refer to \cite{meethal2023finite}. While \cite{meethal2023finite,LeDuc2023} are concerned with numerical parameters, an extension to a still fully discretized, but more operator-learning-like scenario has been presented in \cite{Yamazaki2024} in the context of a time-dependent linear heat equation. Among the vast amount of literature on operator learning \cite{Lu2021}, another popular method to address parametric problems, we mention exemplarily the recent contribution on hypersonic flows \cite{Peyvan2025}.

In our opinion, combining the mathematically proven and quantified ability of FEMs to solve differential equations with the strong approximation capabilities of neural networks to deal with parameter dependencies, should be considered as a simple, but highly appealing idea. In the present work, we thus stick to this approach and extend it in two main directions: 
\begin{enumerate} 
\item First, we consider stationary 2D Stokes and Navier--Stokes equations, i.e., equations of type and characteristics that, given their inherent saddle-point structure, are even for classical FEM quite challenging. 
\item Second, we modify the physics-informed loss function by adding a preconditioner, from left and right, directly to the (nonlinear) residual. In conjunction with L-BFGS, this significantly speeds up training and finally improves accuracy of the obtained solutions. 
\end{enumerate} 

It should be noted that preconditiong from \emph{both} sides and solving the nonlinear problem all-at-once instead of solving a sequence of linear problems distinguishes our approach from the preconditioning strategy proposed in \cite{Songming2024}. Moreover, although the model problem in \cite{Cao2025} is very similar to ours, their ansatz for preconditioning is completely different from ours. 

Regarding the existing literature on physics-informed machine learning for fluid dynamics problems, we refer to the aforementioned survey papers on discretization-free approaches, and only comment on techniques that explicitly rely on discretization. The neural network multigrid solver \cite{Margenberg2022} already mentioned above has been applied to Navier--Stokes equations. DiscretizationNet \cite{Ranade2021} combines finite volume discretization for the instationary Navier--Stokes equations together with an iterative application of an encoder-decoder convolutional neural network, while \cite{Wandel2021} employ physics-informed loss evaluated on a grid together with a 3D U-Net architecture for the neural network. Finally, we mention that beside these machine learning-based approaches, a large number of classical model reduction techniques have been developed over the last decades, especially for applications in fluid dynamics; see, e.g., \cite{BennerMOR}.

\subsubsection*{Structure of the paper}
In Sect.~\ref{sec:method} we introduce the proposed method; first, we describe the underlying setting and the idea in general, and then we discuss this approach in more detail on behalf of our two concrete model problems in Sects.~\ref{subsec:method:stokes} and \ref{subsec:method:navierstokes}. Numerical examples illustrating the advantages of our method are given in Sect.~\ref{sec:numerics}, and a prototypical application to probabilistic inverse problems is presented in Sect.~\ref{sec:application}. We conclude the paper with a brief summary and outlook on future work in Sect.~\ref{sec:conclusion}.  Finally, \ref{sec:proof} contains the proof of a fact that can be understood as a heuristic motivation for the chosen approach.

\section{Preconditioned loss functions for physics-informed, FEM-based neural networks}
\label{sec:method}

In this section, we introduce our proposed method and briefly discuss its realization, both first on an general, rather abstract level. After that, we provide some more concrete details regarding its realization and implementation on the basis of two exemplary, but non-trivial, model problems. 

\subsection{Problem statement} In the following, we will consider problems that can be written in the following abstract variational form: \begin{align}\label{eq::var_form}
    \textnormal{Seek}\; u \in V \qquad \textnormal{s.t.} \quad R(u,v) = 0 \quad \forall v \in V,
\end{align}
where $V$ denotes a suitable vector space\footnote{$V$ is usually a space of functions of a certain regularity, e.g., a Sobolev space.} and $R\!\colon V \times V \rightarrow \R$ is continuous and linear in the second argument, but not necessarily bi- or trilinear. In fact, we will even consider the following parameterized version of \eqref{eq::var_form}: 
\begin{align}\label{eq::var_form_par}
    \textnormal{Seek}\; u(\lambda) \in V \qquad \textnormal{s.t.} \quad R_\lambda(u(\lambda),v) = 0 \quad \forall v \in V,
\end{align}
where $\lambda \in \Lambda \subset \R^p$ denotes a vector of parameters from a parameter range $\Lambda$ and each $R_\lambda\!\colon V \times V \rightarrow \R$ satisfies the assumptions as stated for \eqref{eq::var_form}. As the notation "$u(\lambda)$" indicates, we will \emph{assume} in the following that for each $\lambda \in \Lambda$ there is a unique $u_\lambda \in V$ that satisfies \eqref{eq::var_form_par} which makes the map $\lambda \mapsto u(\lambda)$ a well-defined map $\Lambda \rightarrow V$. 

Let now $V_h \subset V$ be a finite dimensional subspace, e.g., resulting from a conforming finite element discretization of $V$; the index $h$ indicates the discretization parameter, e.g., the mesh size used for finite element discretization. The so-called Galerkin discretization of our problem \eqref{eq::var_form_par} is given by:
\begin{align}\label{eq::var_form_par_discr}
    \textnormal{Seek}\; u_h(\lambda) \in V_h \qquad \textnormal{s.t.} \quad R_\lambda(u_h(\lambda),v_h) = 0 \quad \forall v_h \in V_h.
\end{align}
Again, we will assume without further notice that also this discrete problem is well-defined in the sense that for each $\lambda \in \Lambda$ there is a unique $u_h(\lambda) \in V_h$ satisfying \eqref{eq::var_form_par}. 
To obtain a more familiar notation, one may choose a basis $v_{h}^i \in V_h$, $i = 1,...,N_h$, of $V_h$ such that $u_h(\lambda) = \sum_{i=1}^{N_h} {\bf u}_{h}^i(\lambda) v_{h}^i$ is satisfied with a (parameter-dependent) coefficient vector ${\bf u}_{h}(\lambda) \in \R^{N_h}$. Then, it suffices to test \eqref{eq::var_form_par_discr} with $v_h = v_{h}^i$ for $i=1,...,N_h$ and one obtains a system of $N_h$ equations for ${\bf u}_h(\lambda)$: 
\begin{align}\label{eq::res_form_par_eqn}
\textnormal{Seek}\; {\bf u}_h(\lambda) \in \R^{N_h} \qquad \textnormal{s.t.} \quad {\bf R}_{\lambda,h} ({\bf u}_h(\lambda)) := \left( R_\lambda\big( u_h(\lambda), v_{h}^i \big)\right)_{i=1,...,N_h} = {\bf 0}. 
\end{align}

\subsection{Proposed Method}
Since \eqref{eq::res_form_par_eqn} is a possibly nonlinear system of equations, usually consisting of a high number of variables, solving \eqref{eq::res_form_par_eqn} for several choices of $\lambda$ may become costly. To deal with this issue, we propose a method that makes use of the strong approximation capabilities of neural networks. More precisely, we want to approximate the map $\lambda \mapsto {\bf u}_h(\lambda)$ solving \eqref{eq::res_form_par_eqn} by a neural network ${\bf u}_{h,\NN}^\theta\!\colon \R^p \rightarrow \R^{N_h}$; hereby, $\theta$ denotes the network parameters of the neural network. Training this neural network --- or in other words: finding the right parameters $\theta$ that make ${\bf u}_{h,\NN}^\theta$ a hopefully good approximation of $\lambda \mapsto {\bf u}_h(\lambda)$ --- requires to determine a suitable so-called loss function. Previous work \cite{meethal2023finite} has suggested to use the $\ell^2$-norm of the discrete equation residual, i.e., to determine $\theta = \argmin_{\vartheta} \Lkal$ of 
\begin{equation}\label{eq::eucl_loss}
    \Lkal := \frac{1}{\lvert \Lambda_\train \rvert} \sum_{\lambda \in \Lambda_{\train}} \lVert {\bf R}_{\lambda,h} ({\bf u}_{h,\NN}^\vartheta(\lambda)) \rVert_{\ell^2}^2,
\end{equation}
where $\Lambda_\train \subset \Lambda$ denotes a finite subset used as training data. In this paper we propose a slightly different approach, denoting our NN output as ${\bf \tilde{u}}_{h,\NN}^\vartheta(\lambda)$, we determine $\theta = \argmin_{\vartheta} \Lkal_{\bf P_h}$ via 
\begin{equation}\label{eq::precon_loss}
    \Lkal_{\bf P_h} :=  \frac{1}{\lvert \Lambda_\train \rvert} \sum_{\lambda \in \Lambda_{\train}} \lVert \mathbf{P}^L_{\lambda,h} {\bf R}_{\lambda,h} (\mathbf{P}^R_{\lambda,h} {\bf \tilde{u}}_{h,\NN}^\vartheta(\lambda)) \rVert_{\ell^2}^2,
\end{equation}
with ${\bf P}^L_{\lambda,h}, {\bf P}^R_{\lambda,h} \in \R^{N_h \times N_h}$  being the left, and respectively right, preconditioners. The desired ${\bf u}_{h,\NN}^\theta(\lambda)$ can then be reobtained via ${\bf u}_{h,\NN}^\theta(\lambda) = {\bf P}^R_{\lambda,h} {\bf \tilde{u}}_{h,\NN}^\theta(\lambda)$. Of course, such a preconditioner must be chosen carefully taking into account the problem characteristics. As already mentioned in the introduction, we directly apply preconditioning the the possibly nonlinear residual, which is in contrast to \cite{Songming2024} where nonlinear problems are treated by solving linear problems in iterative manner. Moreover, unlike in \cite{Songming2024} we precondition from both sides.

\begin{figure}[t]
\color{black}\renewcommand\color[2][]{}
\centering
\input{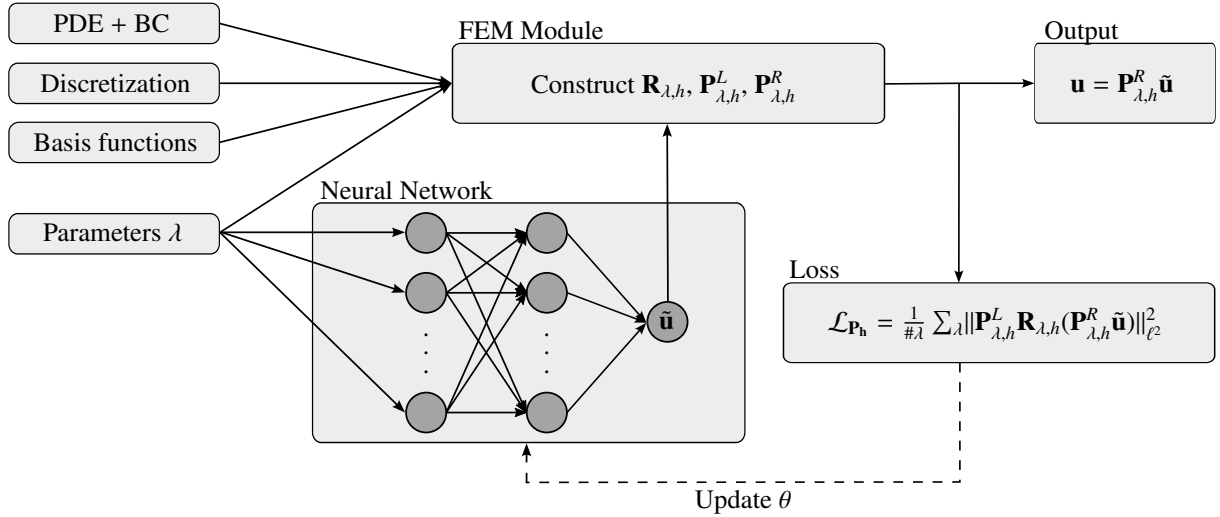}
\caption{Schematic of FEM-based neural networks.}
\label{fenn_schematic}
\end{figure}

Figure~\ref{fenn_schematic} shows a schematic illustration of how FEM-based neural networks can be implemented in practice on a general level. The neural network takes the parameter vector $\lambda$ as input and produces $\tilde{\bf u}_{h,\NN}^\theta(\lambda)$ as output. During training, the FEM module needs to carry out the calculation of the preconditioned residual ${\bf P}^L_{\lambda,h} {\bf R}_{\lambda,h} \circ {\bf P}^R_{\lambda,h}$, either directly and all-at-once, or by calculating/applying the residual and the preconditioners separately; from these quantities the preconditioned loss function $\Lkal_{{\bf P}_{h}}$ is computed. During inference, the FEM module only needs to provide the right-preconditioner ${\bf P}^R_{\lambda,h}$ in order to allow the computation of ${\bf u}_{h,\NN}^\theta(\lambda)$ from the NN output $\tilde{\bf u}_{h,\NN}^\theta(\lambda)$. 

To optimize the loss function w.r.t. the neural network parameters $\theta$, usually variants of the gradient descent method, such as the Adam optimizer \cite{Kingma2014}, or quasi-Newton methods, such as the limited-memory Broyden--Fletcher--Goldfarb--Shanno (L-BFGS) algorithm~\cite{Liu1989} are employed. These methods require the gradient of the loss function with respect to the network parameters, which is calculated by using automatic differentiation (AD)~\cite{griewank2014automatic}; thus all computations in the FEM Module, such as, e.g., setting up stiffness and/or mass matrices etc., must be implemented automatically differentiable. In this paper we will solely employ L-BFGS. Like most optimizers popular in deep learning, L-BFGS only requires gradients of the objective function, i.e., first-order information, and thus avoids the expensive computation of second-order derivatives by AD, which would be necessary for second-order methods as, e.g., Newton's method. Nevertheless, as a quasi-Newton method, L-BFGS lies somewhere between first- and second-order methods as it is able to exploit curvature information by successively building an approximation to the inverse Hessian matrix during the iteration. The associated hopes for faster convergence have been confirmed, e.g., in the case of classic PINNs \cite{RAISSI2019686}.

Our desire of preconditioning the problem is also largely driven by this choice of optimizer. Since~\cite{nazareth1979BFGS} it is known that BFGS for a quadratic loss function and using perfect line search yields the same iterates as a CG algorithm. Since the CG algorithm is known to have a convergence rate that is highly dependent on the condition number of the Hessian, it seems desirable to precondition the optimization problem. In fact, if ${\bf R}_{\lambda,h} (\mathbf{P}^R_{\lambda,h} {\bf \tilde{u}}_{h,\NN}^\theta(\lambda)) = 0$ (or $\approx 0$, respectively) for $\lambda \in \Lambda_\train$ , it holds that 
\[
\nabla^2_\theta \Lkal_{{\bf P}_h} (\theta) = \frac{2}{\lvert \Lambda_\train \rvert} \sum_\lambda H_\lambda^T H_\lambda
\]
(or that $\nabla^2_\theta \Lkal_{{\bf P}_h} (\theta)$ is dominated by this term, respectively) with
\[
H_\lambda = \mathbf{P}^L_{\lambda,h} \cdot [D_{\bf u}{\bf R}_{\lambda,h}] (\mathbf{P}^R_{\lambda,h} {\bf \tilde{u}}_{h,\NN}^\theta(\lambda)) \cdot \mathbf{P}^R_{\lambda,h} \cdot D_\theta {\bf \tilde{u}}_{h,\NN}^\theta(\lambda).
\]
Hence, we hope to improve the conditioning of the Hessian of $\Lkal_{{\bf P}_h}$ at some minimizer $\theta$ by choosing ${\bf P}^L_{\lambda,h}$, ${\bf P}^R_{\lambda,h}$ in such a way that the condition number of $\mathbf{P}^L_{\lambda,h} \cdot [D_{\bf u}{\bf R}_{\lambda,h}] (\mathbf{P}^R_{\lambda,h} {\bf \tilde{u}}_{h,\NN}^\theta(\lambda)) \cdot \mathbf{P}^R_{\lambda,h}$ becomes relatively small. Even if these considerations by no means constitute a strict proof, we consider them to be a conclusive motivation for our approach, and our numerical observations in Sect.~\ref{sec:numerics} seem to support this.

\subsection{Example 1: Stokes equations in 2D} 
\label{subsec:method:stokes}

In the following, we show some more concrete details regarding the realization and implementation of the proposed FEM-based neural network for solving the Stokes flow around an airfoil with a parameterizable angle of attack.
We consider the problem in a 2-dimensional Lipschitz-bounded domain $\Omega$, where the Stokes problem consists of finding two functions, the velocity $u(\lambda)\!: \Omega \rightarrow \mathbb{R} ^2$ and the pressure $p(\lambda)\!: \Omega \rightarrow \mathbb{R}$ such that
\begin{align}
    - \Delta u(\lambda) +  \nabla p(\lambda) &= f_{\lambda},\label{momcons}\\
    \nabla \cdot u(\lambda) &= g_{\lambda},\label{masscons}
\end{align}
with  functions $f_\lambda\!: \Omega \rightarrow \mathbb{R} ^2$, $g_\lambda\!: \Omega \rightarrow \mathbb{R}$ and suitable boundary conditions, which will be introduced precisely below.
Equation~\eqref{momcons} represents the conservation of momentum, Equation~\eqref{masscons} the conservation of mass. Of course, depending on the regularity of the formulation one seeks, these equations can be understood in a strong sense, almost everywhere, or, if $f,g$ become distributions, in a distributional sense.

For the sake of simplicity we will focus in the following on $f\in L^2(\Omega)^2$ and $g\in L^2(\Omega)$, and choose function spaces $X = \set{x \in H^1(\Omega)^2}{\left.x\right|_{\partial\Omega_D} = 0}$ for velocity and $M = L^2(\Omega)$ for pressure. Here, $\partial\Omega_D$ is supposed to signify the Dirichlet boundary, which shall not comprise the entire boundary $\partial\Omega$. For the time being, we impose zero as the Dirichlet boundary condition for the velocity, which will be alleviated later on.

A weak formulation of the system \eqref{momcons}--\eqref{masscons} with a conforming Galerkin ansatz space $V_h := X_h \times M_h \subset X \times M$ then reads as follows: 
\begin{equation}
\left\{
\begin{matrix}
\text{Seek}\; u_h \in X_h\; \text{and}\; p_h \in M_h\; \text{such that}\\
a(u_h,v_h) + b(v_h,p_h) = f(v_h) \quad \forall v_h \in X_h,\\
b(u_h,q_h) = g(q_h) \quad \forall q_h \in M_h,
\end{matrix}
\right.
\label{discrete_s}
\end{equation}
with the bilinear forms $a(u_h,v_h) := \int_\Omega \nabla u_h : \nabla v_h$ and $b(v_h,p_h) := - \int_\Omega p_h \nabla \cdot v_h$. The two functionals $f$ and $g$ are defined as $f(v_h) \coloneqq \int_\Omega v_h f_\lambda$ and $g(q_h) \coloneqq - \int_\Omega q_h g_\lambda$, respectively. $R_{\lambda}$ is accordingly defined as the residual of Equation~\eqref{discrete_s}. Note that the fully continuous variant, that is, setting $X_h = X$ and $M_h = M$, defines a well-posed problem~\cite[Prop.~4.7]{Ern2004theory}. However, to obtain a well-posed problem in the discrete setting, the discrete problem must satisfy the  Ladyzhenskaya--Babu\u{s}ka--Brezzi (LBB) condition \cite{ladyzhenskaya1969mathematical,bab73, brezzi74} and hence, the choice of $X_h, M_h$ must be handled with care. We use Taylor--Hood mixed finite elements, i.e., continuous, piecewise linear ansatz functions for the pressure and continuous, piecewise quadratic ansatz functions for the velocity. With this choice of $V_h$, the LBB condition for \eqref{discrete_s} is satisfied and we obtain a well-posed discrete problem ~\cite[Lemma~4.23]{Ern2004theory}. 

\begin{figure}[t]
\centering
\includegraphics[width=0.3\textwidth]{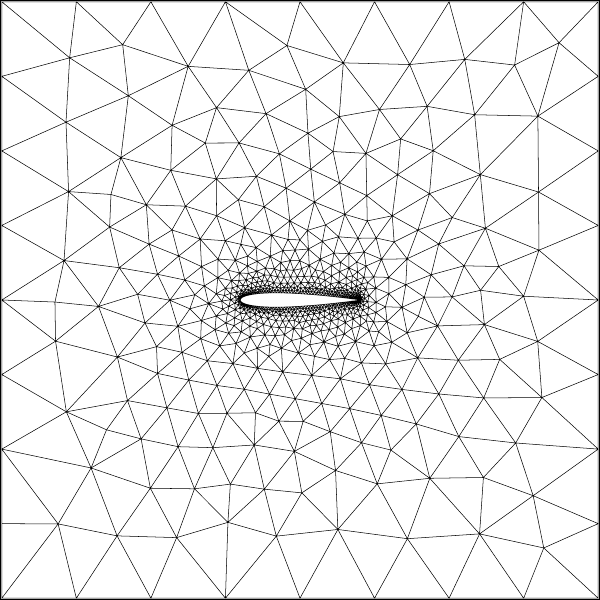}
\caption{Coarse triangular mesh around the NACA 0012 airfoil.}
\label{mesh}
\end{figure}

Let $\{v_h^i\}_{1\leq i \leq N_u}, \{q_h^i\}_{1\leq i \leq N_p}$ be a basis for $X_h, M_h$, with $N_u, N_p$ the respective dimension of the corresponding subspace. When defining $\mathbf{u}= (\mathbf{u}_1,...,\mathbf{u}_{N_u})^T$ with $u_h = \sum_{i=1}^{N_u} \mathbf{u}_i v_h^i \in X_h$ and $\mathbf{p}= (\mathbf{p}_1,...,\mathbf{p}_{N_p})$ with $p_h = \sum_{i=1}^{N_p} \mathbf{p}_i q_h^i \in M_h$, we obtain the residual
\begin{equation}
\mathbf{R}_{\lambda, h} (\mathbf{u},\mathbf{p}) =
    \begin{bmatrix}
 \mathbf{A} & \mathbf{B}^T \\
 \mathbf{B} & \mathbf{0} \\
\end{bmatrix}
\begin{bmatrix}
 \mathbf{u} \\
 \mathbf{p} \\
\end{bmatrix}
-
\begin{bmatrix}
 \mathbf{f} \\
 \mathbf{g} \\
\end{bmatrix},
\label{stokes_mat}
\end{equation}
with $\mathbf{A}_{i,j}= a(v_h^j, v_h^i) \in \mathbb{R}^{N_u, N_u}$, $\mathbf{B}_{k,i}= b(v_h^i, q_h^k) \in \mathbb{R}^{N_p, N_u}$, $\mathbf{f}_i=f(v_h^i) \in \mathbb{R}^{N_u}$ and $\mathbf{g}_k= g(q_h^k) \in \mathbb{R}^{N_p}$. The matrix $\mathbf{A}$ is symmetric positive definite and the system matrix in \eqref{stokes_mat} symmetric, indefinite, and exhibits a typical saddle point structure. For such kind of linear systems, a basic block diagonal preconditioner is given by $\begin{bmatrix}
     \mathbf{A} & \mathbf{0} \\
     \mathbf{0} & -\mathbf{S} \\
\end{bmatrix}^{-1}$, 
where $\mathbf{S} = -\mathbf{BA}^{-1}\mathbf{B}^T$ is the so-called Schur complement; see, e.g., \cite{Benzi2005}. Since $\mathbf{A}$ and $-\mathbf{S}$ are symmetric positive definite, we can use their Cholesky decompositions $\mathbf{A} = \mathbf{LL}^T$ and $-\mathbf{S} = \mathbf{MM}^T$ to define our two preconditioners by 
\begin{align}
    {\bf P}^L_{\lambda,h} = \begin{bmatrix}
     \mathbf{L}^{-1} & \mathbf{0} \\
     \mathbf{0} & \mathbf{M}^{-1} \\
\end{bmatrix}, \qquad {\bf P}^R_{\lambda,h} = \begin{bmatrix}
     \mathbf{L}^{T} & \mathbf{0} \\
     \mathbf{0} & \mathbf{M}^{T} \\
\end{bmatrix}. 
\label{eq:def:precon}
\end{align}
Hence, we obtain the preconditioned residual introduced in Sect.~\ref{sec:method} as follows: 
\begin{align} \nonumber
{\bf P}^L_{\lambda,h} {\bf R}_{\lambda,h} (\mathbf{u},\mathbf{p}) &=
    \begin{bmatrix}
 \mathbf{L}^T & \mathbf{L}^{-1}\mathbf{B}^T \\
 \mathbf{M}^{-1}\mathbf{B} & \mathbf{0} \\
\end{bmatrix}
\begin{bmatrix}
 \mathbf{u} \\
 \mathbf{p} \\
\end{bmatrix}
-
\begin{bmatrix}
 \mathbf{L}^{-1}\mathbf{f} \\
 \mathbf{M}^{-1}\mathbf{g}\\
\end{bmatrix} \\ \nonumber
&=
\begin{bmatrix}
 \mathbf{I} & \mathbf{L}^{-1}\mathbf{B}^T \mathbf{M}^{-T} \\
 \mathbf{M}^{-1}\mathbf{BL}^{-T} & \mathbf{0} \\
\end{bmatrix}
\begin{bmatrix}
 \mathbf{L}^{T}\mathbf{u} \\
 \mathbf{M}^{T}\mathbf{p} \\
\end{bmatrix}
-
\begin{bmatrix}
 \mathbf{L}^{-1}\mathbf{f} \\
 \mathbf{M}^{-1}\mathbf{g}\\
\end{bmatrix}, \\
&=
\begin{bmatrix}
 \mathbf{I} & \mathbf{L}^{-1}\mathbf{B}^T \mathbf{M}^{-T} \\
 \mathbf{M}^{-1}\mathbf{BL}^{-T} & \mathbf{0} \\
\end{bmatrix}
\begin{bmatrix}
 \mathbf{\tilde u} \\
 \mathbf{\tilde p} \\
\end{bmatrix}
-
\begin{bmatrix}
 \mathbf{L}^{-1}\mathbf{f} \\
 \mathbf{M}^{-1}\mathbf{g}\\
\end{bmatrix},
\label{loss_pre}
\end{align}
where $\mathbf I$ is the identity matrix and $\tilde{\bf u} := \mathbf{L}^T {\bf u}$, $\tilde{\bf p} := \mathbf{M}^T {\bf p}$, i.e., $\begin{bmatrix}
 \mathbf{\tilde u} \\
 \mathbf{\tilde p} \\
\end{bmatrix} = {\bf P}_{\lambda,h}^R \begin{bmatrix}
 \mathbf{u} \\
 \mathbf{p} \\
\end{bmatrix}$.  

Lemma~2.1 in \cite{rusten1992preconditioned} shows that the conditioning of the plain saddle point problem \eqref{stokes_mat} is bounded from below by the conditioning of the discretization of the Laplacian. As such, classical FEM theory~\cite[Theorem~9.14]{Ern2004theory} shows that for a Lagrange finite element, as we have chosen for the velocity discretization, the conditioning scales with $h^{-2}$.
Therefore, it is usually necessary to apply preconditioning to the simple saddle-point problem. In fact, the condition number of the square of the preconditioned stiffness matrix in \eqref{loss_pre} is bounded independently of the discretization, as shown in~\ref{sec:proof}. There, it is even shown that the number of distinct eigenvalues of this matrix is three. Considering again that BFGS with perfect line search is equivalent to CG for quadratic optimization problems~\cite{nazareth1979BFGS}, BFGS applied to such an quadratic optimization problem with the squared preconditioned matrix as its Hessian, would converge in at most three iterations\footnote{The situation is analogous to the Schur-preconditioned stable FEM solutions of the Stokes problem, when solved with Krylov subspace methods, such as MINRES, cf.~\cite[p.~292]{Elman2005}.}. Of course, these properties do not directly transfer to our setting, as the problem under consideration is highly non-linear and non-convex due to the presence of a neural network. Nevertheless, it does not seem unreasonable to assume that applying the preconditioner to the linear part of the problem will lead to at least a gradual improvement of the overall properties in this case as well. In fact, our numerical examples in Sect.~\ref{sec:numerics} support this.

Let us now turn to the explicit numerical example, where we want to solve the Stokes equations \eqref{momcons}-\eqref{masscons} with source terms $f_\lambda = 0, g_\lambda =0$, and non-homogeneous Dirichlet boundary conditions. We define the domain $\Omega$ as  $(0,5)^2$ with a NACA 0012 airfoil \cite{Jacobs1932TheCO} of length one cut out in the center  (see Figure~\ref{mesh}).
A Dirichlet boundary condition dependent on the angle of attack $\lambda$ is used at the inflow boundary $\partial\Omega^+ := \{(x,y)\;\vert\; x =0\; \text{or}\; y=0\}$, a Neumann boundary condition at the outflow boundary $\partial\Omega^- := \{(x,y)\;\vert\; x =5\; \text{or}\; y=5 \}$ and a no-slip boundary condition at the airfoil boundary $\partial\Omega^A$:
\begin{align*}
    u &= [\cos(\lambda \pi/180),\textnormal{ } \sin(\lambda \pi/180)]^T &\;\text{on}\; \partial\Omega^+,\\
    p - \partial_\nu u &= 0 &\;\text{on}\; \partial\Omega^-,\\
    u &= 0 &\;\text{on}\; \partial\Omega^A.
\end{align*}
Here $\partial_\nu$ signifies the partial derivative in normal direction on the boundary.
Note that the homogeneous Neumann boundary condition on the outflow boundary effectively sets the mean outflow pressure to zero.

Figure~\ref{mesh} shows the triangular mesh used in the experiments. As we expect larger changes closer to the airfoil, the mesh closer to the airfoil is finer than at the outer edge. We use $2\times3193$ degrees of freedom for the velocity and 839 degrees of freedom for the pressure, which corresponds to 7225 degrees of freedom in total. 

Homogeneous Neumann boundary conditions are known to be so-called natural conditions and thus included in~\eqref{stokes_mat} by construction, and do not need any special treatment. However, our aforementioned theory assumes homogeneous Dirichlet boundary conditions for the velocity as well. In order to enforce the inhomogeneous boundary conditions, we will implicitly construct a lifting function: Firstly, we assemble all the matrices and residuals oblivious to the Dirichlet boundary conditions. Secondly, we exploit that the chosen basis functions are Lagrangian, i.e., each basis function $v_h^i$ has an associated Lagrangian node $x_i \in \partial\Omega$ such that
\begin{equation}
    v_h^i(x_j)=\delta_{ij} \quad \forall i,j.
    \label{lagrange}
\end{equation}
The Dirichlet boundary conditions are now enforced on all such $x_j\in \partial\Omega^+ \cup \partial\Omega^{A}$ by adjusting the right-hand side vector $\mathbf f$ and the stiffness matrix in~\eqref{stokes_mat} accordingly: First, we multiply the columns of the stiffness matrix corresponding to the boundary nodes by the specified boundary values and then subtract the result from $\mathbf f$. Next, we set the corresponding values of $\mathbf f$ to the boundary values. In the stiffness matrix we modify the rows and columns corresponding to the boundary nodes so that the diagonal value is unity and the off-diagonal entries are set to zero. %
Note that this is done before preconditioning, which is possible because the modifications do not affect the symmetry properties of the system~\eqref{stokes_mat}.

In the present case, the (preconditioned) stiffness matrix and right hand side vector of the residual \eqref{stokes_mat} are independent of the prediction of the neural network. Thus they need to be calculated only once and can be reused for every loss calculation in the training of the neural network.

It is also noteworthy that, by the aforementioned approach of incorporating the Dirichlet boundary conditions, we obtain in principle a single objective, avoiding the issues of weakly imposing the Dirichlet conditions via an additional objective term. However, we believe that avoiding the multi-objective scenario is not per se an issue, as one could also rewrite our objective with multiple terms. In our view, the more important part is that these potentially additional terms that impose the Dirichlet conditions must be mutually non-conflicting. We explicitly ensure this by erasing parts of the assembled stiffness matrix and adequately modifying the source terms, as mentioned above.

\subsection{Example 2: Stationary Navier--Stokes equations in 2D}
\label{subsec:method:navierstokes}

\begin{figure}[t]
\centering
    \renewcommand\sffamily{}
    \input{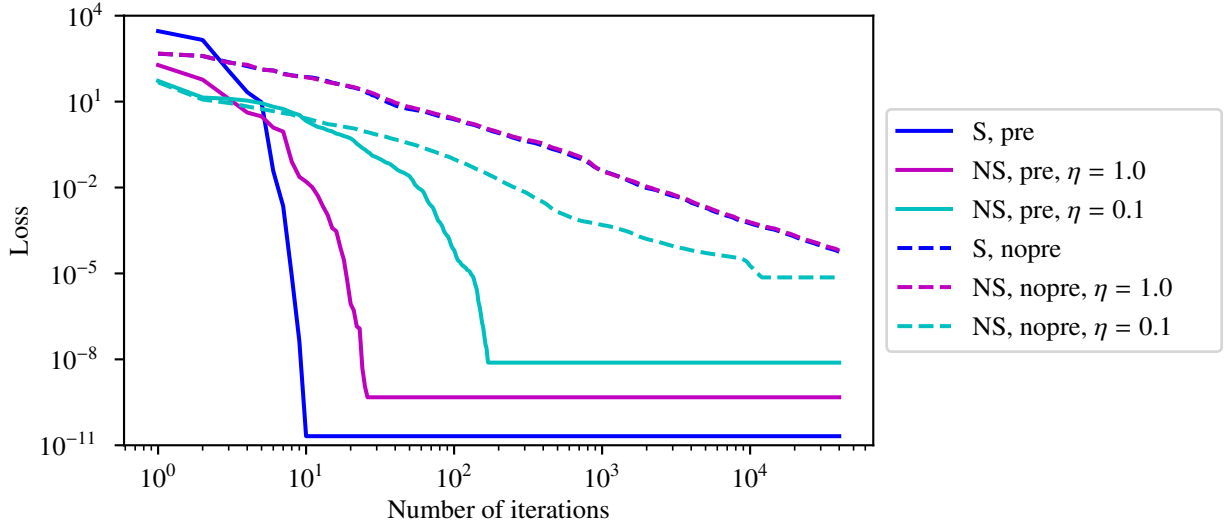}
\caption{Loss values per L-BFGS iteration when training FEM-based NNs with preconditioning (pre) and without preconditioning (nopre) for Stokes (S) and Navier--Stokes (NS) equations with an angle of attack $\lambda = 1^{\circ}$ and different viscosities $\eta$.
\label{fig:loss_O=1}}
\end{figure}

To avoid redundancies, we focus on the changes that are necessary to extend our approach from the Stokes to the Navier--Stokes equations. Thus, details not discussed further remain the same as in the previous section. The Navier--Stokes equations are as follows: 

\begin{align}
    - \eta \Delta u(\lambda) + u(\lambda) \cdot \nabla u(\lambda) + \nabla p(\lambda) &= f_{\lambda},\label{ns_mom}\\
    \nabla \cdot u(\lambda) &= g_{\lambda},\label{ns_mass}
\end{align}
i.e., compared to the Stokes equations \eqref{momcons}-\eqref{masscons} in particular the non-linear convection term $u(\lambda) \cdot \nabla u(\lambda)$ has been added. Moreover, we have introduced the so-called kinematic viscosity $\eta>0$ as additional parameter. Discretization leads to 
\begin{equation}
\left\{
\begin{matrix}
\text{Seek}\; u_h \in X_h\; \text{and}\; p_h \in M_h\; \text{such that}\\
a(u_h,v_h) + c(z_h;u_h,v_h) + b(v_h,p_h) = f(v_h) \quad \forall v_h \in X_h,\\
b(u_h,q_h) = g(q_h) \quad \forall q_h \in M_h,
\end{matrix}
\right.
\label{discrete_ns}
\end{equation}
with $b$ as before, $a$ as before but multiplied by $\eta$, and the trilinear form $c(z_h; u_h, v_h) = \int_\Omega (z_h \cdot \nabla u_h) \cdot v_h$ originating from the nonlinear term. We again use Taylor--Hood elements for the discretization in the non-linear case, as these are already necessary for a stable solution in the linear case.
To formulate the corresponding nonlinear algebraic equation, we utilize Lagrangian basis functions as in \eqref{stokes_mat} and obtain the discrete residual
\begin{equation}
\mathbf{R}_{\lambda, h} (\mathbf{u},\mathbf{p}) =
    \begin{bmatrix}
 \mathbf{A} + \mathbf{C}(\mathbf{u}) & \mathbf{B}^T \\
 \mathbf{B} & \mathbf{0} \\
\end{bmatrix}
\begin{bmatrix}
 \mathbf{u} \\
 \mathbf{p} \\
\end{bmatrix}
-
\begin{bmatrix}
 \mathbf{f} \\
 \mathbf{g} \\
\end{bmatrix},
\label{ns_mat}
\end{equation}
with $\mathbf{C}(\mathbf{u})_{i,j}= c(u_h;v_h^j, v_h^i) \in \mathbb{R}^{N_u, N_u}$.
We use the same preconditioning as for the Stokes problem, which yields the preconditioned system
\begin{equation}
{\bf P}_{\lambda,h} {\bf R}_{\lambda,h} (\mathbf{u},\mathbf{p}) =
    \begin{bmatrix}
 \mathbf{I} +  \mathbf{L}^{-1}\mathbf{C}(\mathbf{u}) \mathbf{L}^{-T} & \mathbf{L}^{-1}\mathbf{B}^T \mathbf{M}^{-T} \\
 \mathbf{M}^{-1}\mathbf{BL}^{-T} & \mathbf{0} \\
\end{bmatrix}
\begin{bmatrix}
 \mathbf{L}^{T}\mathbf{u} \\
 \mathbf{M}^{T}\mathbf{p} \\
\end{bmatrix}
-
\begin{bmatrix}
 \mathbf{L}^{-1}\mathbf{f} \\
 \mathbf{M}^{-1}\mathbf{g}\\
\end{bmatrix}.
\label{loss_pre_ns}
\end{equation}
Again we use the NN to directly predict $\tilde{\bf u} := \mathbf{L}^T {\bf u}$ and $\tilde{\bf p} := \mathbf{M}^T {\bf p}$. In our numerical examples we will solve the Navier--Stokes equations \eqref{ns_mom}-\eqref{ns_mass} on the same geometry, mesh, and boundary conditions as in the previous section. 

Since the Navier--Stokes problem is nonlinear, the stiffness matrix depends on the solution, i.e., the prediction of the NN in our case. This means that unlike in \eqref{loss_pre}, where the system matrix and the right hand side could be calculated once in advance, we now have to recalculate parts of the residual \eqref{loss_pre_ns} in every L-BFGS iteration during the training of the network. More precisely, this applies to the convection term $\mathbf{C}(\mathbf{u})$ and the right hand side vector $\mathbf{f}$ due to the inclusion of Dirichlet boundary conditions as we will describe below. We note that
\begin{align*}
\mathbf{C}_{ij}(\mathbf{u}) = \int_\Omega \sum_k \mathbf{u}_k v_h^k \cdot \nabla v_h^j \cdot v_h^i = \sum_k \mathbf{u}_k \underbrace{\int_\Omega  v_h^k \cdot \nabla v_h^j \cdot v_h^i}_{=: \mathbf{\Tilde{C}}_{i,j,k}}, 
\end{align*}
i.e., $\mathbf{C}(\mathbf{u})$ can be calculated by a contraction of the tensor $\mathbf{u} = \mathbf{L}^{-T} \tilde{\bf u}$ (or $\mathbf{u} = \tilde{\bf u}$ in the case of no preconditioning) with a 3D-tensor $\mathbf{\Tilde{C}}$ that is independent of $\mathbf{u}$ and thus can be precomputed and reused. Since $\mathbf{\Tilde{C}}$ is sparse, this can be done in a memory-efficient manner.

Let us now briefly comment on enforcing Dirichlet boundary conditions in the residual \eqref{ns_mat}. Except for the nonlinear part, $\mathbf{C}(\mathbf{u})$, this can be handled in the same way as for the Stokes equations. The only difference is that instead of $\mathbf{C}(\mathbf{u})$ we actually have to compute $\mathbf{C}(\hat{\mathbf{u}})$ with the slightly modified argument $\hat {\bf u} = \hat {\bf u}({\bf u})$ in which the respective boundary conditions are already enforced at all Lagrange nodes on the Dirichlet boundary, i.e.,
\begin{equation}
    \mathbf{\hat{u}}_i := \left\{ 
    \begin{array}{ll} 
    \mathbf{u}_i & \text{if}\; x_i \notin \partial\Omega^+\cup \partial\Omega^A, \\
    u_{\textnormal{BC}}(x_i) & \text{if}\; x_i \in \partial\Omega^+\cup \partial\Omega^A,
    \end{array}
    \right.    
\end{equation}
where $u_{\textnormal{BC}}(x_i)$ denotes the value of $u$ prescribed at the boundary node $x_i$. 
As for the Stokes equations, the boundary conditions are included into \eqref{ns_mat} before preconditioning. With these changes, the training of the FEM-based NN can be carried out as described above.

\section{Numerical results}
\label{sec:numerics}

In the following, we show results of the proposed FEM-based neural networks for solving the Stokes and Navier--Stokes flow around an airfoil with a parameterizable angle of attack as described in Section \ref{sec:method}. In particular, we will focus on the comparison between using a preconditioner and not using a preconditioner.

The respective neural networks and an FEM module as described in Sects.~\ref{sec:method}, \ref{subsec:method:stokes}, and \ref{subsec:method:navierstokes} have been implemented in PyTorch \cite{PyTorch2019}. Unless stated otherwise, the neural networks are plain vanilla feed forward networks consisting of five hidden layers with 50 neurons each and an output layer whose size is equal to the number of degrees of freedom of the discretization. When testing various hyperparameters, we found that small changes in the number of layers and number of neurons per layer did not have a major effect on the results. However, it turned out that a non-linear activation function which is not bounded is of great advantage to avoid local minima during optimization, which is why we use SELU~\cite{klambauer2017self} as the activation function except for the last, linear layer. To facilitate efficient training, we set the network weights with Xavier's normal initialization \cite{glorot2010understanding} scaled by $3/4$.
For the training, we use PyTorch's implementation\footnote{\href{https://pytorch.org/docs/stable/generated/torch.optim.LBFGS.html}{https://pytorch.org/docs/stable/generated/torch.optim.LBFGS.html} [Accessed August 21, 2024]} of L-BFGS. To make the results more comparable with ``classical" implementations, we choose {\tt max\_iter=1} to ensure that only one L-BFGS step is performed per iteration. Furthermore, we choose a line search within L-BFGS that seeks a solution that fulfills the strong Wolfe condition.

In the following, the term ``error" does not refer to the discretization error, but to the error between our NN-based solution and a pure FEM-based reference solution that uses the same mesh and Taylor--Hood discretization. The latter was computed with FEniCSx \cite{BarattaEtal2023, AlnaesEtal2014}. 

\begin{figure}[t]
\centering
    \renewcommand\sffamily{}
    \input{eigenvalues_hess.pgf}
\caption{Eigenvalues of the Hessian of the loss function, evaluated after training the FEM-based NNs with preconditioning (pre) and without preconditioning (nopre) for Stokes equations with an angle of attack $\lambda = 1^{\circ}$. Eigenvalues are shown in gray, a kernel density estimate of the respective distribution functions in blue and orange. Note that the y-axis employs symmetric log axis scaling, i.e. logarithmic scaling with a linear scaling exception for values within $[-2,2]$.
\label{fig:eigenvalues}}
\end{figure}

\subsection{The non-parametric case}
We start with experiments with only a single angle of attack as training data and use a linear NN consisting of only one input and one output layer. Considering this non-parametric and thus highly underdetermined case allows us to evaluate the benefits of preconditioning separately from other effects. 

Both the preconditioned and the non-preconditioned NNs are trained and evaluated for a single angle of attack $\lambda = 1^{\circ}$ and 40000 L-BFGS iterations, for each: the Stokes problem and the Navier--Stokes equation. In the latter case, the solutions are computed with viscosities $\eta=1$ and $\eta=0.1$. Figure~\ref{fig:loss_O=1} shows that in all these cases preconditioning significantly speeds up the training in terms of a drastically reduced number of iterations. For the preconditioned cases we can even observe a stagnation of the residual at levels that are probably related to numerical accuracy; this is probably related to a similar reason as the effect described in the paragraph below \eqref{loss_pre} in Sect.~\ref{subsec:method:stokes}. Figure~\ref{fig:time_O=1} shows for the Navier--Stokes equation with $\eta=1.0$ that training is also much faster in terms of runtime because the very significant reduction in the number of iterations more than offsets the minimal increase in the runtime of a single iteration that is due to the few additional operations related to preconditioning.

\begin{figure}[t]
\centering
    \renewcommand\sffamily{}
    \input{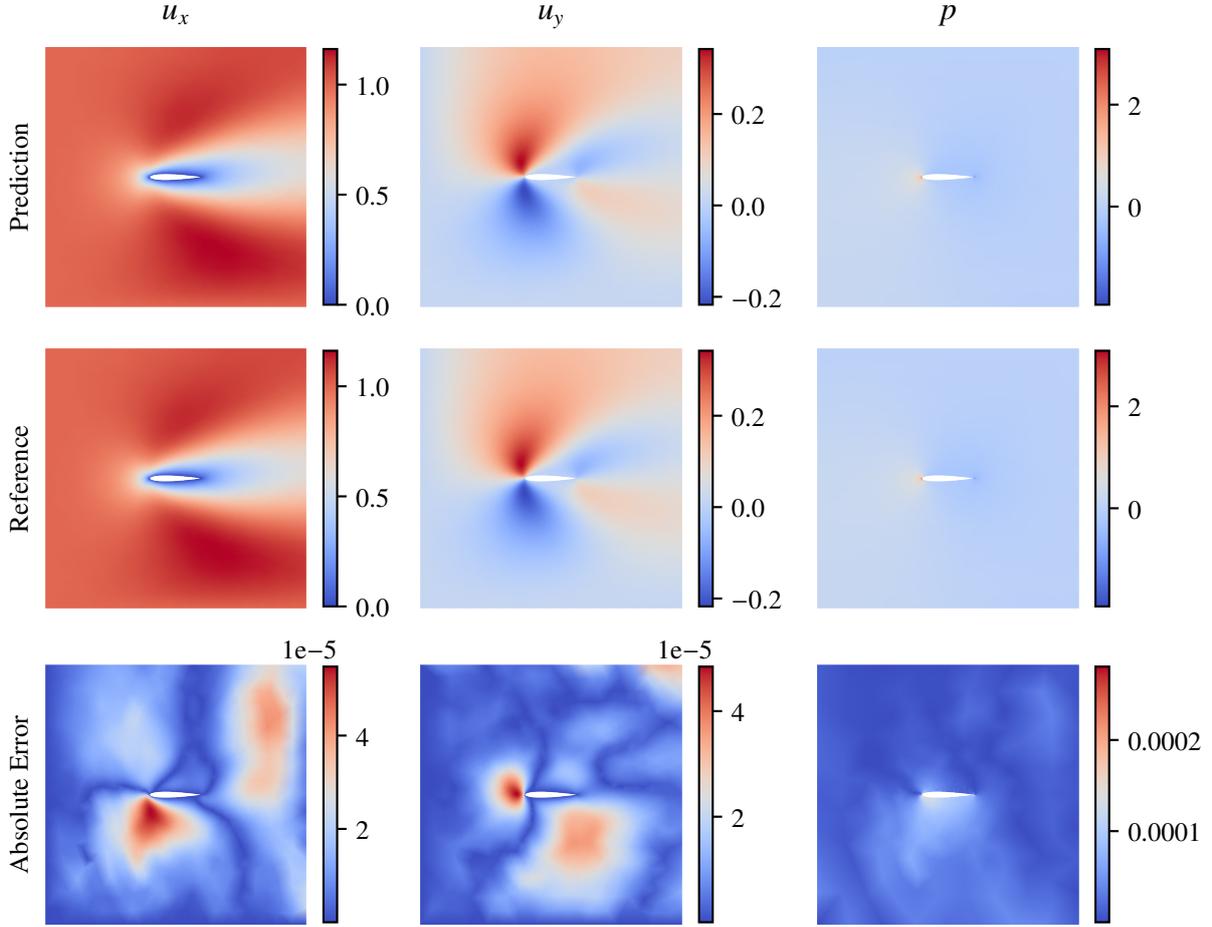}
\caption{Results for Navier--Stokes equation trained with an angle of attack $\lambda = 1^{\circ}$ and viscosity $\eta = 0.1$. The first row shows the prediction of the FEM-based neural network with preconditioning, the second row the reference solution obtained with the FEM, the third row the absolute error of the prediction.}
\label{fenn_results_ns_pre}
\end{figure}

Table~\ref{tab:errors} shows the achieved relative errors in $L^2$- and $L^\infty$-norms. Again, the benefits of preconditioning are clearly visible in the reduction of the final error by several orders of magnitude. In particular, note that for $\eta=0.1$ (which corresponds to an increased Reynolds number) a satisfying accuracy could only be achieved with preconditioning. As Figure~\ref{fenn_results_ns_pre} illustrates, there is almost no visible difference between the NN-based solutions and the pure FEM-based ones in the end.

To further illustrate the effectiveness of the preconditioning, we computed the eigenvalues of the Hessian at the last iteration of the training of the Stokes equations, once with preconditioning and once without preconditioning. Note that computing the Hessian using the automatic differentiation machinery in PyTorch is possible, though possibly quite costly, and can only be employed for small problem sizes. The eigenvalues, and a kernel density estimate, for the two cases, with and without preconditioning, is plotted in Figure~\ref{fig:eigenvalues}. It is clearly visible that even though our argumentation in Section~\ref{sec:method} and \ref{sec:proof} totally neglected the non-linearity introduced by the neural networks, the conclusion of clustered eigenvalues still holds. As noted earlier, we assume L-BFGS coupled with perfect line search to share the same characteristics as the CG method for a nearly quadratic optimization problem. Given that in our least-squares setting the cluster of zero eigenvalues can be neglected~\cite{kammerer1972convergence}, the convergence of the CG algorithm for a quadratic optimization problem is determined by the number of remaining eigenvalue clusters~\cite{Elman2005}. Of course, this analysis is not exact, given that our problem is made non-linear by the neural network and potentially the advection term in the Navier--Stokes equations, and thus potentially deviates a lot from the aforementioned assumptions. Nonetheless, it clearly demonstrates the mechanisms at work, indicating the effectiveness of classical preconditioning techniques.

That the assumptions made in the design of the preconditioner are only approximately fulfilled, could not be more visible than in the high Reynolds number case, as illustrated in Figure~\ref{fig:loss_O=1}. For $\eta = 0.01$ the advection term clearly dominates, such that the chosen preconditioner, which takes only the diffusive part into account, performs poorly. The latter is clearly indicated by the increased number of iterations. Moreover, for higher Reynolds numbers techniques like, e.g, Streamline Upwind Petrov--Galerkin (SUPG) should be employed to further enhance the numerical properties~\cite{Elman2005,Brooks1982}.

\begin{figure}[t]
\centering
    \renewcommand\sffamily{}
    \input{time_comparison_O=1.pgf}
\caption{Loss values over time, in logarithmic scaling, when training FEM-based NNs with preconditioning (pre) and without preconditioning (nopre) for the Navier--Stokes (NS) equations with an angle of attack $\lambda = 1^{\circ}$ and viscosity $\eta = 1.0$. The initial construction times of the system matrices in \eqref{ns_mat} and \eqref{loss_pre_ns} are shown in cyan, the training times of the iterations in purple.
\label{fig:time_O=1}}
\end{figure}

\begin{table}[t]
    \centering
    \begin{tabular}{ |c|c|c|r r r|r r r| }
 \hline
    & & & \multicolumn{3}{|c|}{$L^{2}$-Error} & \multicolumn{3}{|c|}{$L^{\infty}$-Error}\\
  & Pre. & $\eta$ & \multicolumn{1}{c}{$u_x$} & \multicolumn{1}{c}{$u_y$} & \multicolumn{1}{c|}{$p$} & \multicolumn{1}{c}{$u_x$} & \multicolumn{1}{c}{$u_y$} & \multicolumn{1}{c|}{$p$} \\ 
 \hline
 \multirow{2}{*}{S} & nopre & - & 6.0801e-03 & 2.5899e-02 & 8.6753e-02 & 7.6672e-02 & 1.1719e-01 & 9.9965e-01 \\

 & pre & - & 7.3656e-08 & 1.1760e-07 & 1.3784e-07 & 1.5618e-07 & 1.6593e-07 & 4.1410e-06 \\ 
 \hline
\multirow{4}{*}{NS} & nopre & 1.0 & 7.3695e-03 & 4.4370e-02 & 1.5622e-01 & 8.3582e-02 & 1.5751e-01 & 9.9945e-01 \\
 & nopre & 0.1 & 8.8126e-01 & 15.1430e-00 & 24.3703e-00 & 1.7591e-00 & 11.3491e-00 & 2.2046e-00 \\

 & pre & 1.0 & 1.4740e-06 & 8.7947e-06 & 5.9979e-06 & 3.1987e-06 & 9.7846e-06 & 7.3386e-06 \\ 

  & pre & 0.1 & 1.7492e-05 & 1.2692e-04 & 2.0684e-04 & 4.7161e-05 & 1.4152e-04 & 9.0761e-05 \\ 
 \hline
 \end{tabular}
 \phantom{1}
\caption{Relative errors of the FEM-based neural network trained for Stokes (S) and Navier--Stokes (NS) using the the non-preconditioned loss \eqref{eq::eucl_loss} (nopre) and the preconditioned loss \eqref{eq::precon_loss} (pre) with an angle of attack of 1° and viscosity $\eta$.}
    \label{tab:errors}
\end{table}

\subsection{The parametric case}

In this section, we consider the true parametric problem and train with multiple angles of attack, namely 3, 5, 9, and 17 different $\lambda$ distributed equidistantly in the range $[1^{\circ},45^{\circ}]$. A maximum of 4000 iterations of L-BFGS is used to train the respective networks. 

\begin{figure}[t]

    \renewcommand\sffamily{}
    \centering
    \input{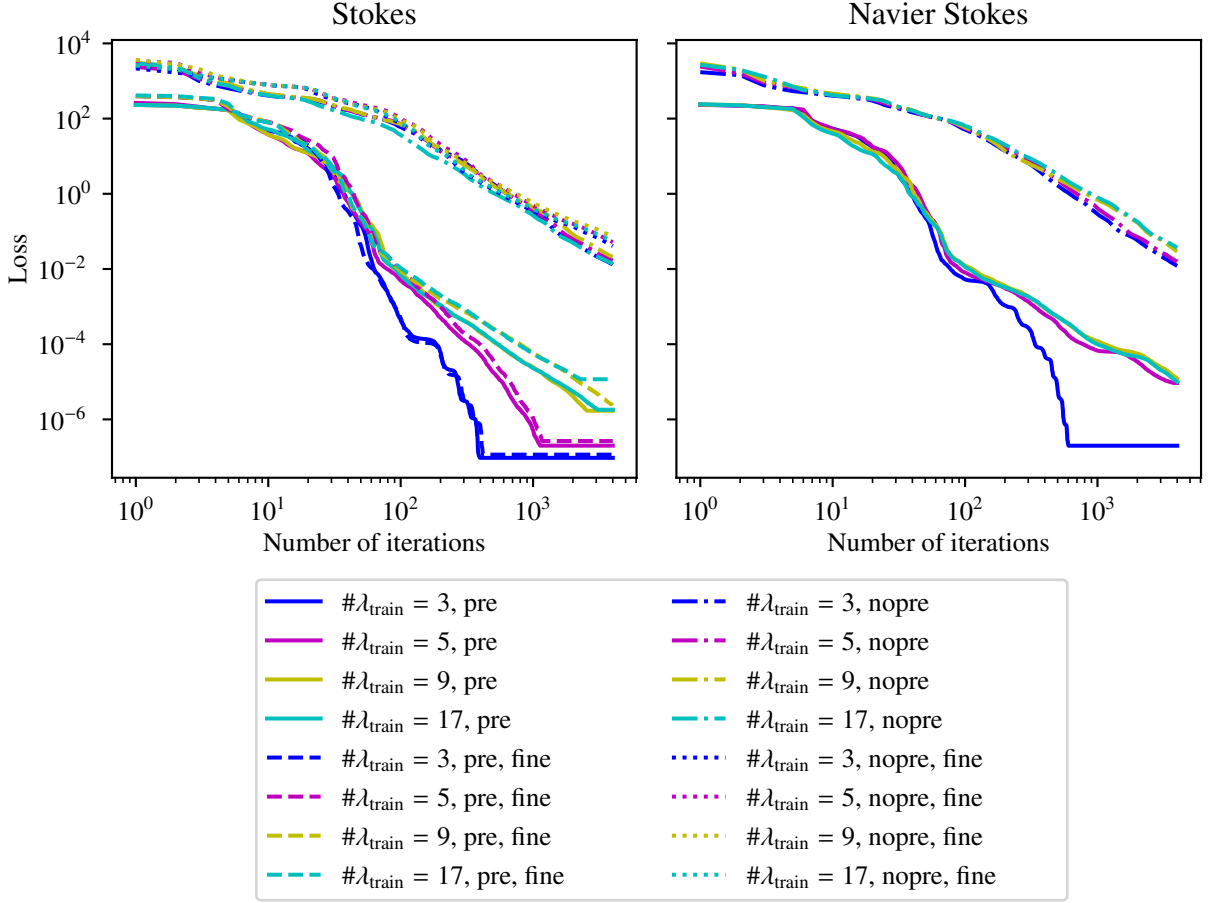}
    \caption{Loss values per L-BFGS iteration when training FEM-based NNs for multiple angles of attack $\#\lambda_\train$ with preconditioning (pre) and without preconditioning (nopre). For the Stokes problem, in addition to the loss behavior for the original mesh, we also show the behavior for a finer mesh.}
    \label{loss_stokes_O>0}
\end{figure}

Figure~\ref{loss_stokes_O>0} shows the loss values during training for the different configurations under consideration. Even in this case, preconditioning is able to accelerate the decay of the loss by more than two orders of magnitude for both the Stokes and the Navier--Stokes equations. For the Stokes problem, we also performed the experiments on a finer mesh, which is obtained from the original one by a single uniform refinement and results in 28085 degrees of freedom. The behavior of the losses on the original and the finer mesh can be taken as an indication of mesh independence, which is consistent with our expectations from~\ref{sec:proof} and our discussion in the previous section. 

\begin{figure}[t]
    \renewcommand\sffamily{}
    \centering \input{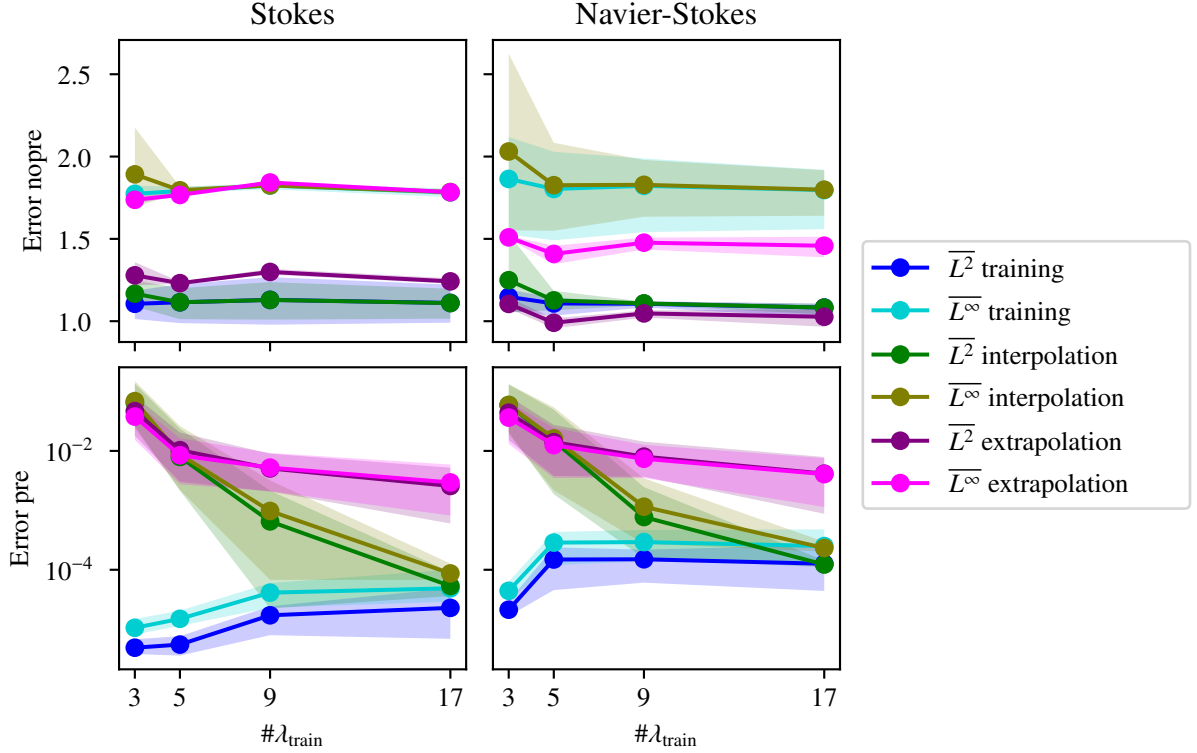}
    \caption{Averaged relative $L^2$- and $L^{\infty}$-errors of the variable $u_x$ predicted from FEM-based NNs trained  with and without preconditioning for Stokes and Navier--Stokes ($\eta = 1$) with several numbers of angles of attack $\#\lambda_{\text{train}}$ and tested with angles of attack inside (interpolation) and outside (extrapolation) the training range. The dots represent the mean errors, the shadings correspond to the range in which the errors occurred for the given set of testing angles.}
    \label{fenn_generalizabilty}
\end{figure}

In the following, we will in particular consider the training errors, i.e., the error on the seen training data, as well as generalization errors: more specifically the errors seen on interpolating input data as well as on extrapolating input data. To evaluate the interpolation and extrapolation qualities, we test with unseen angles of attack within the training range, $\lambda \in \{5^{\circ}, 16.5^{\circ}, 30^{\circ}, 40^{\circ}\}$, and outside the training range, $\lambda \in \{47.5^{\circ}, 50^{\circ}, 55^{\circ}\}$, respectively. 
Figure~\ref{fenn_generalizabilty} shows the average errors on the training data ("training") as well as the generalization errors ("interpolation" and "extrapolation") in $L^2$- and $L^\infty$-norm for the different sizes of the training set. For clarity, we only report errors of the velocity in $x$-direction as the errors of velocity in $y$-direction and pressure behave similarly. 
It is clearly visible that preconditioning does not only heavily improve quality of the surrogate model on the training set itself, but also its inter- and extrapolation capabilities; in fact, the errors obtained without preconditioning are far beyond an acceptable range. The increase of the training errors for a growing size of the training set for the preconditioned problems is likely due to the fact that the same network architecture (i.e. the same number of degrees of freedom) is used to solve a increasingly complex problem. As expected, the generalization errors decrease with increasing training data.

Our findings can be summarized as follows: The use of preconditioned loss functions speeds up the training procedure, both for linear and nonlinear problems and, moreover, improves the quality of the resulting surrogate models both in terms of training and generalization errors. 

\section{Application: Solving an inverse problem related to stationary 2D Navier--Stokes equations}
\label{sec:application}

In the following, we demonstrate the relevance of our method by applying it to a prototypical inverse problem related to the stationary 2D Navier--Stokes equations. This problem corresponds to an application scenario in which pressure measurements from two sensors on the airfoil are used to determine the angle of attack. Since measurements in a real-world setting are always noisy, it is crucial to be able to quantify the aleatoric uncertainty coming from this noise. Our goal is to estimate a probability density for the unknown attack angle $\lambda$ associated with these observations. Therefore, we pursue a probabilistic approach and use a Markov Chain Monte Carlo (MCMC) method in combination with our fully differentiable FEM-based NNs.

\begin{figure}[t]
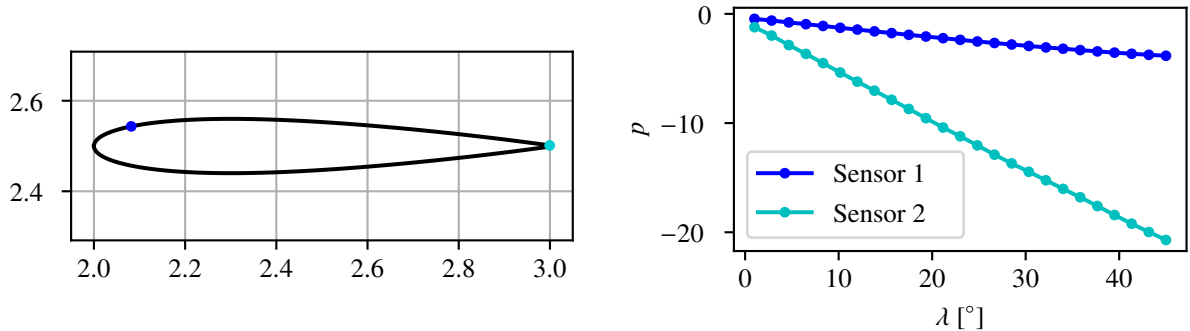


    \renewcommand\sffamily{}
    \begin{subfigure}[c]{0.49\textwidth}
    \centering
    \input{Figure7_1.pgf}
    \end{subfigure}
    \begin{subfigure}[c]{0.49\textwidth}
    \input{Figure7_2.pgf}
    \end{subfigure}
    \caption{Left: Location of sensors on airfoil. Right: Pressure $p$ at sensor locations for different angles of attack $\lambda$ predicted by FEM-based NN.}
    \label{fig:p_meas}
\end{figure}

Let $\mathbf{y}_1, \ldots, \mathbf{y}_m \in \R^2$ be noisy observations of the pressure at two fixed points of the airfoil. We assume that the two pressure sensors are located at 8.18\% and 99.94\% of the chord on the top of the airfoil; see Figure~\ref{fig:p_meas} for an illustration. Moreover, we assume that the measurements are independent and randomly distributed according to a normal distribution with standard deviation $\sigma = 0.5$. We thus obtain a log-likelihood that modulo some $\lambda$-independent constant $C$ is given by
\begin{equation}
    \log p(y|\lambda) = C - \sum_{i=1}^m \frac{1}{\sigma^2} \lVert \mathbf{\hat{p}}_\lambda - \mathbf{y}_i\rVert^2\, .
\end{equation}
We deliberately use $\mathbf{\hat{p}}_\lambda$ to signify the result of the preconditioned neural network evaluated at just these two measurement points. Encompassing this likelihood with a uniform prior $p(\lambda) \sim \mathcal{U}([1^{\circ},45^{\circ}])$, we are interested in the resulting posterior distribution $p(\lambda|y) \propto p(y|\lambda) p(\lambda)$.

Given that computations of $\mathbf{\hat{p}}_\lambda$ exhibit a nonlinear dependence on $\lambda$, we employ a Hamiltonian Monte Carlo (HMC) method \cite{betancourt2013, neal2012}, more precisely the No-U-Turn Sampler (NUTS). The latter is a gradient-based MCMC method that approximates the desired density function. HMC allows to sample from this posterior by using an approximate Hamiltonian dynamics simulation, which is then corrected by a Metropolis acceptance step. To solve the Hamiltonian differential equation, derivatives of the target posterior density function with respect to $\lambda$ must be computed. Here, an important advantage of neural network-based surrogate modeling comes into play: the surrogate model is not only differentiable (which may theoretically also be the case with other techniques), but in fact these derivatives can be computed efficiently using AD. We use the implementation of NUTS that comes with the PyTorch-based library Pyro~\cite{bingham2018pyro}.

\begin{figure}[t]
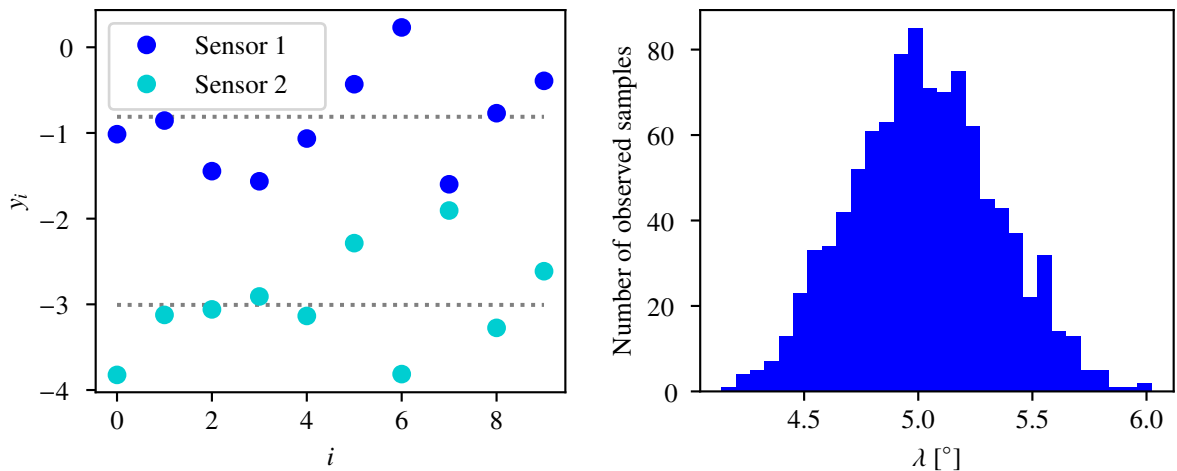


    \renewcommand\sffamily{}
    \centering
    \input{Figure8_1.pgf}
    \input{Figure8_2.pgf}
    \caption{Left: Observed pressure measurements $y_i$ of two sensors. Right: Corresponding posterior distribution of the angle of attack $\lambda$ of the airfoil, as generated with the MCMC algorithm.}
    \label{uq}
\end{figure}

In order to test this method, we choose the following setup: Figure~\ref{fig:p_meas} (right) shows the forward problem of predicting the pressure at the sensor locations for 17 angles of attack $\lambda$ that are equidistantly distributed in the range $[1^{\circ},45^{\circ}]$.
Using this FEM-based NN as surrogate model, we generate 10 artificial and noisy pressure ``measurements" by adding Gaussian noise with a standard deviation of $0.5$ to the respective pressure predictions evaluated at a reference angle of attack of $\lambda = 5^{\circ}$. The ``measurements" are shown in Figure~\ref{uq}.

Figure~\ref{uq} displays the 1000 samples of the posterior distribution, generated by MCMC with 1000 warm-up steps. The generated distribution has an empirical mean of 5.03 and a standard deviation of 0.32. Hence, the combination of MCMC and FEM-based NNs was able to deliver plausible results for this test case.

\section{Conclusion}
\label{sec:conclusion}

FEM-based/FEM-inspired neural networks are known to combine the best of two worlds: finite elements and neural networks. On the one hand, the (parametric) PDE under consideration is discretized by a method that is well-understood and reliable, the incorporation of boundary conditions and complex geometries is relatively straightforward. On the other hand, the discretization benefits from the power of neural networks to handle parameter dependence, e.g., by obtaining a cheap-to-evaluate and fully differentiable parameter-to-solution map, which can be advantageous in the context of inverse problems. 

In this paper, we have advanced this type of technique by applying it for the first time to PDEs from fluid dynamics that result in a saddle point problem after discretization. In particular, we have introduced a modification of this approach in which preconditioning applied to the residual loss allows to improve both the speed of training and the resulting accuracy. We have demonstrated the relevance of our approach with numerical examples on behalf of a prototypical but realistic model problem from fluid dynamics and a related inverse problem. 

Having provided a proof-of-concept for our method in this paper, an interesting goal for future research would be to investigate the possibilities for the inclusion of, e.g., distributed parallel methods or low-cost preconditioners, as they are typical for large-scale finite element implementations. Another goal would be to consider other PDEs or PDE systems with different characteristics (including suitable preconditioners). In particular, the extension to time-dependent problems or to flow problems with higher Reynolds numbers that require stabilization shall be mentioned.

\bigskip

\appendix 
\section{Mesh-independent conditioning of the FEM layer for the preconditioned Stokes case}
\label{sec:proof}

In this appendix, we show that contribution of the matrix in~\eqref{loss_pre} to the loss in~\eqref{eq::precon_loss} is well-conditioned, as far as the last layer of the NN is concerned. More precisely, we will show that the square of
\begin{align}\label{eqn:defz}
    \mathbf{Z} &= \begin{bmatrix}
 \mathbf{I} & \mathbf{L}^{-1}\mathbf{B}^T \mathbf{M}^{-T} \\
 \mathbf{M}^{-1}\mathbf{BL}^{-T} & \mathbf{0} \\
\end{bmatrix}
\end{align}
has an eigenvalue distribution that is independent of the mesh resolution. The following derivation is based on~\cite{Murphy2000}.

\begin{lem}
Consider $\mathbf{Y} = \mathbf{Z}^2$ with $\mathbf{Z}$ as given in \eqref{eqn:defz}, then $\mathbf{Y}$ has the following annihilating polynomial
\begin{gather*}
    \left(\mathbf{Y} - \mathbf{I}\right)\left(\left(\mathbf{Y}-\frac{3}{2}\mathbf{I}\right)^2 - \frac{5}{4}\mathbf{I}\right) = 0\, .
\end{gather*}
\end{lem}
\begin{proof}
Introducing $\mathbf{C} = \mathbf{M}^{-1}\mathbf{B}\mathbf{L}^{-T}$, first note that $\mathbf{Y}$ is given by
\begin{align*}
    \mathbf{Y} &= \begin{bmatrix}
 \mathbf{I} + \mathbf{C}^T \mathbf{C} & \mathbf{C}^T \\
 \mathbf{C} & \mathbf{C}\, \mathbf{C}^T \\
\end{bmatrix} \, .
\end{align*}
The latter can, using the definitions of $\mathbf{L},\mathbf{M},\mathbf{S}$, be simplified to
\begin{align*}
    \mathbf{C}\, \mathbf{C}^T =& \mathbf{M}^{-1} \mathbf{B} \underbrace{\mathbf{L}^{-T} \mathbf{L}^{-1}}_{=\mathbf{A}^{-1}} \mathbf{B}^T \mathbf{M}^{-T}\\
        =& \mathbf{M}^{-1} \underbrace{\mathbf{B} \mathbf{A}^{-1} \mathbf{B}^T}_{-\mathbf{S}} \mathbf{M}^{-T}\\
        =& \mathbf{I}\, .
\end{align*}
Furthermore, $\mathbf{P} = \mathbf{C}^T \mathbf{C}$ is a projection operator
\begin{align*}
    \mathbf{C}\, \mathbf{P} =& \mathbf{C}\, \mathbf{C}^T \mathbf{C} = \mathbf{C}\\
    \mathbf{P}^2 =& \mathbf{C}^T \mathbf{C}\, \mathbf{P} = \mathbf{C}^T \mathbf{C} = \mathbf{P}\, .
\end{align*}
As such, we obtain
\begin{align*}
    \left(\mathbf{Y} -\frac{3}{2}\mathbf{I}\right)^2 =& \begin{bmatrix}
                                         \left(\mathbf{P}-\frac{1}{2}\mathbf{I}\right)^2 + \mathbf{P} & 0 \\
                                         0 & \frac{5}{4}\mathbf{I} \\
                                        \end{bmatrix}\, ,
\end{align*}
which yields
\begin{align*}
    \left(\mathbf{Y} -\frac{3}{2}\mathbf{I}\right)^2 - \frac{5}{4}\mathbf{I} =& \begin{bmatrix}
                                         \mathbf{P}-\mathbf{I} & 0 \\
                                         0 & 0 \\
                                        \end{bmatrix}\, .
\end{align*}
The conclusion then follows easily from the observation that
\begin{align*}
    \mathbf{Y} - \mathbf{I} &= \begin{bmatrix}
 \mathbf{P} & \mathbf{C}^T \\
 \mathbf{C} & 0 \\
\end{bmatrix} \, .
\end{align*}
\end{proof}

As an immediate consequence of the last Lemma, the distinct eigenvalues of $\mathbf{Y}$ are $\left\{1, \frac{3}{2} \pm \frac{\sqrt{5}}{2}\right\}$, thus clearly bounding the condition number of the matrix $\mathbf{Y}$.

\section*{Acknowledgements}
This research was carried out in part during the project PISA (Physics Inspired AI) by the German Aerospace Center (DLR). The DLR  is supported by the Federal Ministry for Economic Affairs and Climate Action (BMWK) on the basis of a decision by the German Bundestag.

Last but not least, the authors want to thank the anonymous reviewers for their valuable comments.

\bibliographystyle{unsrt}
\bibliography{references}

\end{document}